# Interval Valued Bipolar Fuzzy Weighted Neutrosophic Sets and Their Application

*Irfan Deli*
Muallim Rıfat Faculty of Education,
Kilis 7 Aralık University,
79000 Kilis, Turkey,
irfandeli@kilis.edu.tr

*Yusuf Şubaş*
Muallim Rıfat Faculty of Education,
Kilis 7 Aralık University,
79000 Kilis, Turkey,
ysubas@kilis.edu.tr

*Florentin Smarandache*
University of New Mexico, 705 Gurley Ave., Gallup,
New Mexico 87301, USA.
fsmarandache@gmail.com

*Mumtaz Ali*
Department of Mathematics,
Quaid-i-Azam University
Islamabad, Pakistan
mumtazali7288@gmail.com

*Abstract*—Interval valued bipolar fuzzy weighted neutrosophic set(IVBFWN-set) is a new generalization of fuzzy set, bipolar fuzzy set, neutrosophic set and bipolar neutrosophic set so that it can handle uncertain information more flexibly in the process of decision making. Therefore, in this paper, we propose concept of IVBFWN-set and its operations. Also we give the IVBFWN-set average operator and IVBFWN-set geometric operator to aggregate the IVBFWN-sets, which can be considered as the generalizations of some existing ones under fuzzy, neutrosophic environments and so on. Finally, a decision making algorithm under IVBFWN environment is given based on the given aggregation operators and a real example is used to demonstrate the effectiveness of the method.

## I. INTRODUCTION

To overcome containing various kinds of uncertainty, the concept of fuzzy sets [18] has been introduced by Zadeh. After Zadeh, many studies on mathematical modeling have been developed. For example; to model indeterminate and inconsistent information Smarandache [13] introduced the concept of neutrosophic set which is independently characterized by three functions called truth-membership function, indeterminacy-membership function and falsity membership function. Recently, studies on neutrosophic sets are made rapidly in [1,2].

Bipolar fuzzy sets, which are a generalization of Zadeh's fuzzy sets [18], were originally proposed by Lee [9]. Bosc and Pivert [4] said that "Bipolarity refers to the propensity of the human mind to reason and make decisions on the basis of positive and negative effects. Positive information states what is possible, satisfactory, permitted, desired, or considered as being acceptable. On the other hand, negative statements express what is impossible, rejected, or forbidden. Negative preferences correspond to constraints, since they specify which values or objects have to be rejected (i.e., those that do not satisfy the constraints), while positive preferences correspond to wishes, as they specify which objects are more desirable than others (i.e., satisfy user wishes) without rejecting those that do not meet the wishes." Presently, works on bipolar fuzzy sets are progressing rapidly in [3,4,8-12,17]. Also, bipolar neutrosophic set(BN-set) and its operations is given in [7].

In this study, to handling some uncertainties in fuzzy sets and neutrosophic sets, the extensions of fuzzy sets[18], bipolar fuzzy sets[9], neutrosophic sets[13] and bipolar neutrosophic sets[7], interval valued bipolar fuzzy weighted neutrosophic sets with application are introduced.

## II. PRELIMINARIES

In the section, we give some concepts related to bipolar fuzzy sets, neutrosophic sets, interval valued neutrosophic set, and bipolar neutrosophic sets.

**Definition 2.1.** [14] Let $X$ be a universe of discourse. Then a single valued neutrosophic set is defined as:

$$A_{NS} = \{\langle x, T_A(x), I_A(x), F_A(x)\rangle : x \in X\}$$

which is characterized by a truth-membership function $T_A(x): X \to [0,1]$, an indeterminacy-membership function $I_A(x): X \to [0,1]$, and a falsity-membership function $F_A(x): X \to [0,1]$. There is not restriction on the sum of $T_A(x), I_A(x),$ and $F_A(x)$ so $0 \le T_A(x) \le I_A(x) \le F_A(x) \le 3$.

**Definition 2.2.** [15] Let $X$, be a space of points (objects) with generic elements in $X$, denoted by $x$. An interval valued neutrosophic set (for short IVNS) $A$ in $X$, is characterized by truth-membership function $T_A(x)$, indeteminacy-membership function $I_A(x)$, and falsity-membership function $F_A(x)$. For each point $x$ in $X$, we have that $T_A(x), I_A(x), F_A(x) \subseteq [0,1]$.

For two IVNS
$$A_{IVNS} = \{\langle x, [\inf T_A(x), \sup T_A(x)], [\inf I_A(x), \sup I_A(x)], [\inf F_A(x), \sup F_A(x)]\rangle : x \in X\}$$
and
$$B_{IVNS} = \{\langle x, [\inf T_B(x), \sup T_B(x)], [\inf I_B(x), \sup I_B(x)], [\inf F_B(x), \sup F_B(x)]\rangle : x \in X\}$$

Then,
1. $A_{IVNS} \subseteq B_{IVNS}$ if and only if
   $\inf T_A(x) \leq \inf T_B(x)$, $\sup T_A(x) \leq \sup T_B(x)$,
   $\inf I_A(x) \geq \inf I_B(x)$, $\sup I_A(x) \geq \sup I_B(x)$,
   $\sup F_A(x) \geq \sup F_B(x)$, $\sup F_A(x) \geq \sup F_B(x)$
   for all $x \in X$.

2. $A_{IVNS} = B_{IVNS}$ if and only if
   $\inf T_A(x) = \inf T_B(x)$, $\sup T_A(x) = \sup T_B(x)$,
   $\inf I_A(x) = \inf I_B(x)$, $\sup I_A(x) = \sup I_B(x)$,
   $\sup F_A(x) = \sup F_B(x)$, $\sup F_A(x) = \sup F_B(x)$
   for any $x \in X$.

3. $A_{IVNS}{}^C$ if and only if
   $A_{IVNS}{}^C = \{\langle x, [\inf F_A(x), \sup F_A(x)], [1 - \sup I_A(x), 1 - \inf I_A(x)], [\inf T_A(x), \sup T_A(x)]\rangle : x \in X\}$

4. $A_{IVNS} \cap B_{IVNS}$ if and only if
   $A_{IVNS} \cap B_{IVNS} = \{\langle x, [\inf T_A(x) \wedge \inf T_B(x), \sup T_A(x) \wedge \sup T_B(x)], [\inf I_A(x) \vee \inf I_B(x), \sup I_A(x) \vee \sup I_B(x)], [\inf F_A(x) \vee \inf F_B(x), \sup F_A(x) \vee \sup F_B(x)]\rangle : x \in X\}$

5. $A_{IVNS} \cup B_{IVNS}$ if and only if
   $A_{IVNS} \cup B_{IVNS} = \{\langle x, [\inf T_A(x) \vee \inf T_B(x), \sup T_A(x) \vee \sup T_B(x)], [\inf I_A(x) \wedge \inf I_B(x), \sup I_A(x) \wedge \sup I_B(x)], [\inf F_A(x) \wedge \inf F_B(x), \sup F_A(x) \wedge \sup F_B(x)]\rangle : x \in X\}$

**Definition 2.3**. [9] Let $X$ be a non-empty set. Then, a bipolar-valued fuzzy set, denoted by $A_{BF}$ is defined as;
$$A_{BF} = \{\langle x, \mu_B^+(x), \mu_B^-(x)\rangle : x \in X\}$$

Where $\mu_B^+(x) : X \to [0,1]$ and $\mu_B^-(x) : X \to [0,1]$. The positive membership degree $\mu_B^+(x)$ denotes the satisfaction degree of an element $x$ to the property corresponding to $A_{BF}$ and the negative membership degree $\mu_B^-(x)$ denotes the satisfaction degree of $x$ to some implicit counter property of $A_{BF}$.

**Definition 2.4.** [7] A bipolar neutrosophic set $A$ in $X$ is defined as an object of the form
$$A = \{\langle x, T^+(x), I^+(x), F^+(x), T^-(x), I^-(x), F^-(x)\rangle : x \in X\},$$
where $T^+, I^+, F^+ : X \to [1,0]$ and $T^-, I^-, F^- : X \to [-1,0]$.

The positive membership degree $T^+(x), I^+(x), F^+(x)$ denotes the truth membership, indeterminate membership and false membership of an element $x \in X$ corresponding to a bipolar neutrosophic set $A$ and the negative membership degree $T^-(x), I^-(x), F^-(x)$ denotes the truth membership, indeterminate membership and false membership of an element $x \in X$ to some implicit counter-property corresponding to a bipolar neutrosophic set $A$.

**Definition 2.5.** [7] Let
$$A_1 = \{\langle x, T_1^+(x), I_1^+(x), F_1^+(x), T_1^-(x), I_1^-(x), F_1^-(x)\rangle : x \in X\}$$
and
$$A_2 = \{\langle x, T_2^+(x), I_2^+(x), F_2^+(x), T_2^-(x), I_2^-(x), F_2^-(x)\rangle : x \in X\}$$
be two bipolar neutrosophic sets.

I. Then $A_1 \subseteq A_2$ if and only if
   $T_1^+(x) \leq T_2^+(x)$, $I_1^+(x) \leq I_2^+(x)$, $F_1^+(x) \geq F_2^+(x)$,
   and
   $T_1^-(x) \geq T_2^-(x)$, $I_1^-(x) \geq I_2^-(x)$, $F_1^-(x) \leq F_2^-(x)$
   for all $x \in X$.

II. Then $A_1 = A_2$ if and only if
    $T_1^+(x) = T_2^+(x)$, $I_1^+(x) = I_2^+(x)$, $F_1^+(x) = F_2^+(x)$
    and
    $T_1^-(x) = T_2^-(x)$, $I_1^-(x) = I_2^-(x)$, $F_1^-(x) = F_2^-(x)$
    for all $x \in X$.

III. Then their union is defined as:
     $(A_1 \cup A_2)(x) =$
     $\langle x, \{\max\{T_1^+(x), T_2^+(x)\}, \frac{I_1^+(x) + I_2^+(x)}{2}, \min\{F_1^+(x), F_2^+(x)\},$
     $\min\{T_1^-(x), T_2^-(x)\}, \frac{I_1^-(x) + I_2^-(x)}{2}, \max\{F_1^-(x), F_2^-(x)\} \} \rangle : x \in X\}$

     for all $x \in X$.

IV. Then their intersection is defined as:

$(A_1 \cap A_2)(x) =$

$\{\langle x, \min\{T_1^+(x), T_2^+(x)\}, \frac{I_1^+(x)+I_2^+(x)}{2}, \max\{F_1^+(x), F_2^+(x)\},$

$\max\{T_1^-(x), T_2^-(x)\}, \frac{I_1^-(x)+I_2^-(x)}{2}, \min\{F_1^-(x), F_2^-(x)\} \rangle : x \in X\}$

for all $x \in X$.

V. Then the complement of $A_1$ is denoted by $A_1^c$ and is defined by

$T_{A_1^c}^+(x) = \{1^+\} - T_{A_1}^+(x), \quad I_{A_1^c}^+(x) = \{1^+\} - I_{A_1}^+(x),$

$F_{A_1^c}^+(x) = \{1^+\} - F_{A_1}^+(x)$

and

$T_{A_1^c}^-(x) = \{1^-\} - T_{A_1}^-(x), \quad I_{A_1^c}^-(x) = \{1^-\} - I_{A_1}^-(x),$

$F_{A_1^c}^-(x) = \{1^-\} - F_{A_1}^-(x),$

for all $x \in X$.

**Definition 2.6.** [7] Let

$A_1 = \langle T_1^+, I_1^+, F_1^+, T_1^-, I_1^-, F_1^- \rangle$

and

$A_2 = \langle T_2^+, I_2^+, F_2^+, T_2^-, I_2^-, F_2^- \rangle$

be two bipolar neutrosophic number. Then the operations for these numbers are defined as below;

a. $\lambda A_1 = \langle 1-(1-T_1^+)^\lambda, (I_1^+)^\lambda, (F_1^+)^\lambda, -(-T_1^-)^\lambda, -(-I_1^-)^\lambda,$

$-(1-(1-(-F_1^-))^\lambda) \rangle$

b. $A_1^\lambda = \langle (T_1^+)^\lambda, 1-(1-I_1^+)^\lambda, 1-(1-F_1^+)^\lambda,$

$-(1-(1-(-T_1^-))^\lambda), -(-I_1^-)^\lambda, -(-F_1^-)^\lambda \rangle$

c. $A_1 + A_2 = \langle T_1^+ + T_2^+ - T_1^+.T_2^+, I_1^+ I_2^+, F_1^+ F_2^+, -T_1^- .T_2^-,$

$-(-I_1^- - I_2^- - I_1^- .I_2^-), -(-F_1^- - F_2^- - F_{1L}^- .F_2^-) \rangle$

d. $A_1 A_2 = \langle T_1^+ T_2^+, I_1^+ + I_2^+ - I_1^+ .I_2^+, F_1^+ + F_2^+ - F_1^+ .F_2^+,$

$-(-T_1^- - T_2^- - T_1^- .T_2^-), I_1^- I_2^-, F_1^- F_2^- \rangle$

where $\lambda > 0$.

**Definition 2.7.** [7] Let

$\tilde{a} = \langle T^+, I^+, F^+, T^-, I^-, F^- \rangle$

be a bipolar neutrosophic number. Then, the score function $s(\tilde{a})$, accuracy function $a(\tilde{a})$ and certainty function $c(\tilde{a})$ of an NBN are defined as follows:

$s(\tilde{a}) = \frac{1}{6}(T^+ + 1 - I^+ + 1 - F^+ + 1 + T^- - I^- - F^-)$

$a(\tilde{a}) = T^+ - F^+ + T^- - F^-$

$c(\tilde{a}) = T^+ - F^-$

**Definition 2.8.** [7] Let

$\tilde{a}_j = \langle T_j^+, I_j^+, F_j^+, T_j^-, I_j^-, F_j^- \rangle (j = 1, 2, ..., n)$

be a family of bipolar neutrosophic numbers. Then,

a) $F_W : \mathfrak{I}_n \to \mathfrak{I}$ is called bipolar neutrosophic weighted average operator if it satisfies;

$F_W(\tilde{a}_1, \tilde{a}_2, ..., \tilde{a}_n) = \sum_{j=1}^{n} w_j \tilde{a}_j$

$= \langle 1 - \prod_{j=1}^{n}(1-T_j^+)^{w_j}, \prod_{j=1}^{n}(I_j^+)^{w_j}, \prod_{j=1}^{n}(F_j^+)^{w_j}, -\prod_{j=1}^{n}(-T_j^-)^{w_j},$

$-\left(1 - \prod_{j=1}^{n}(1-(-I_j^-))^{w_j}\right), -\left(1 - \prod_{j=1}^{n}(1-(-F_j^-))^{w_j}\right) \rangle$

where $w_j$ is the weight of $\tilde{a}_j (j = 1, 2, ..., n)$, $w_j \in [0,1]$ and

$\sum_{j=1}^{n} w_j = 1$.

b) $H_W : \mathfrak{I}_n \to \mathfrak{I}$ is called bipolar neutrosophic weighted geometric operator if it satisfies;

$H_W(\tilde{a}_1, \tilde{a}_2, ..., \tilde{a}_n) = \prod_{j=1}^{n} \tilde{a}_j^{w_j}$

$= \langle \left[\prod_{j=1}^{n}(T_j^+)^{w_j}, 1 - \prod_{j=1}^{n}(1-I_j^+)^{w_j}, 1 - \prod_{j=1}^{n}(1-F_j^+)^{w_j},\right.$

$-\left(1 - \prod_{j=1}^{n}(1-(-T_j^-))^{w_j}\right), -\prod_{j=1}^{n}(-I_{jL}^-)^{w_j}, -\prod_{j=1}^{n}(-F_j^-)^{w_j} \rangle$

where $w_j$ is the weight of $\tilde{a}_j (j = 1, 2, ..., n)$, $w_j \in [0,1]$ and

$\sum_{j=1}^{n} w_j = 1$.

III. INTERVAL VALUED BIPOLAR FUZZY WEIGHTED NEUTROSOPHIC SET

In this section we give concept of IVBFWN-set and its operations. Also we give the IVBFWN-set average operator and IVBFWN-set geometric operator with properties to aggregate the IVBFWN-sets based on the study given in [7].

**Definition 3.1.** A interval valued bipolar fuzzy weighted neutrosophic set(IVBFWN-set) $A$ in $X$ is defined as an object of the form

$A = \{\langle x, [T_L^+(x), T_R^+(x)], [I_L^+(x), I_R^+(x)], [F_L^+(x), F_R^+(x)],$
$[T_L^-(x), T_R^-(x)], [I_L^-(x), I_R^-(x)], [F_L^-(x), F_R^-(x)], p(x)\rangle : x \in X\}$

where $T_L^+, T_R^+, I_L^+, I_R^+, F_L^+, F_R^+ : X \to [0,1]$ and $T_L^-, T_R^-, I_L^-, I_R^-, F_L^-, F_R^- : X \to [-1,0]$. Also $p : X \to [0,1]$ fuzzy weighted index of the element $x$ in $X$.

**Example 3.2.** Let $X = \{x_1, x_2, x_3\}$. Then

$$A = \begin{cases} \langle x_1, [0.3, 0.9], [0.1, 0.8], [0.2, 0.5], [-0.8, -0.7], [-0.5, -0.1], [-0.4, -0.3], 0.5\rangle, \\ \langle x_2, [0.3, 0.8], [0.3, 0.9], [0.1, 0.2], [-0.7, -0.6], [-0.6, -0.2], [-0.6, -0.2], 0.7\rangle, \\ \langle x_3, [0.4, 0.7], [0.5, 0.7], [0.3, 0.4], [-0.9, -0.5], [-0.4, -0.3], [-0.8, -0.1], 0.8\rangle \end{cases}$$

is a IVBFWN subset of $X$.

**Theorem 3.3.** A IVBFWN-set is the generalization of a bipolar fuzzy set and bipolar neutrosophic set.

Proof: Straightforward.

**Definition 3.4.** Let
$A_1 = \{\langle x, [T_{1L}^+(x), T_{1R}^+(x)], [I_{1L}^+(x), I_{1R}^+(x)], [F_{1L}^+(x), F_{1R}^+(x)],$
$[T_{1L}^-(x), T_{1R}^-(x)], [I_{1L}^-(x), I_{1R}^-(x)], [F_{1L}^-(x), F_{1R}^-(x)], p_1(x)\rangle : x \in X\}$

and

$A_2 = \{\langle x, [T_{2L}^+(x), T_{2R}^+(x)], [I_{2L}^+(x), I_{2R}^+(x)], [F_{2L}^+(x), F_{2R}^+(x)],$
$[T_{2L}^-(x), T_{2R}^-(x)], [I_{2L}^-(x), I_{2R}^-(x)], [F_{2L}^-(x), F_{2R}^-(x)], p_2(x)\rangle : x \in X\}$
be two IVBFWN-sets.

1. Then $A_1 \subseteq A_2$ if and only if

$T_{1L}^+(x) \le T_{2L}^+(x), \quad T_{1R}^+(x) \le T_{2R}^+(x), \quad I_{1L}^+(x) \ge I_{2L}^+(x),$
$I_{1R}^+(x) \ge I_{2R}^+(x), \quad F_{1L}^+(x) \ge F_{2L}^+(x), \quad F_{1R}^+(x) \ge F_{2R}^+(x),$
$T_{1L}^-(x) \le T_{2L}^-(x), \quad T_{1R}^-(x) \le T_{2R}^-(x), \quad I_{1L}^-(x) \ge I_{2L}^-(x),$
$I_{1R}^-(x) \ge I_{2R}^-(x), \quad F_{1L}^-(x) \ge F_{2L}^-(x), \quad F_{1R}^-(x) \ge F_{2R}^-(x),$
and
$p_1(x) \le p_2(x)$
for all $x \in X$.

2. Then $A_1 = A_2$ if and only if

$T_{1L}^+(x) = T_{2L}^+(x), \quad T_{1R}^+(x) = T_{2R}^+(x), \quad I_{1L}^+(x) = I_{2L}^+(x),$
$I_{1R}^+(x) = I_{2R}^+(x), \quad F_{1L}^+(x) = F_{2L}^+(x), \quad F_{1R}^+(x) = F_{2R}^+(x),$
$T_{1L}^-(x) = T_{2L}^-(x), \quad T_{1R}^-(x) = T_{2R}^-(x), \quad I_{1L}^-(x) = I_{2L}^-(x),$
$I_{1R}^-(x) = I_{2R}^-(x), \quad F_{1L}^-(x) = F_{2L}^-(x), \quad F_{1R}^-(x) = F_{2R}^-(x),$
and $p_1(x) = p_2(x)$ for all $x \in X$.

3. Then their union is defined as:

$(A_1 \cup A_2)(x) =$
$\{\langle x, [\max\{T_{1L}^+(x), T_{2L}^+(x)\}, \max\{T_{1R}^+(x), T_{2R}^+(x)\}],$
$\left[\frac{I_{1L}^+(x) + I_{2L}^+(x)}{2}, \frac{I_{1R}^+(x) + I_{2R}^+(x)}{2}\right],$
$[\min\{F_{1L}^+(x), F_{2L}^+(x)\}, \min\{F_{1R}^+(x), F_{2R}^+(x)\}],$
$[\min\{T_{1L}^-(x), T_{2L}^-(x)\}, \min\{T_{1R}^-(x), T_{2R}^-(x)\}],$
$\left[\frac{I_{1L}^-(x) + I_{2L}^-(x)}{2}, \frac{I_{1R}^-(x) + I_{2R}^-(x)}{2}\right],$
$[\max\{F_{1L}^-(x), F_{2L}^-(x)\}, \max\{F_{1R}^-(x), F_{2R}^-(x)\}],$
$\max\{p_1(x), p_2(x)\}\rangle : x \in X\}$

for all $x \in X$.

4. Then their intersection is defined as:

$(A_1 \cap A_2)(x) =$
$\{\langle x, [\min\{T_{1L}^+(x), T_{2L}^+(x)\}, \min\{T_{1R}^+(x), T_{2R}^+(x)\}],$
$\left[\frac{I_{1L}^+(x) + I_{2L}^+(x)}{2}, \frac{I_{1R}^+(x) + I_{2R}^+(x)}{2}\right],$
$[\max\{F_{1L}^+(x), F_{2L}^+(x)\}, \max\{F_{1R}^+(x), F_{2R}^+(x)\}],$
$[\max\{T_{1L}^-(x), T_{2L}^-(x)\}, \max\{T_{1R}^-(x), T_{2R}^-(x)\}],$
$\left[\frac{I_{1L}^-(x) + I_{2L}^-(x)}{2}, \frac{I_{1R}^-(x) + I_{2R}^-(x)}{2}\right],$
$[\min\{F_{1L}^-(x), F_{2L}^-(x)\}, \min\{F_{1R}^-(x), F_{2R}^-(x)\}],$
$\min\{p_1(x), p_2(x)\}\rangle : x \in X\}$
for all $x \in X$.

5. Then the complement of $A_1$ is denoted by $A_1^c$, is defined by

$A_1^C = \{\langle x, [F_L^+(x), F_R^+(x)], [1 - I_R^+(x), 1 - I_L^+(x)], [T_L^+(x), T_R^+(x)],$
$[F_L^-(x), F_R^-(x)], [1 - I_R^-(x), 1 - I_L^-(x)], [T_R^-(x), T_L^-(x)], 1 - p(x)\rangle : x \in X\}$

**Example 3.5.** Let $X = \{x_1, x_2, x_3\}$. Then

$$A_1 = \begin{cases} \langle x_1, [0.3, 0.9], [0.1, 0.8], [0.2, 0.5], [-0.8, -0.7], [-0.5, -0.1], [-0.4, -0.3], 0.5\rangle, \\ \langle x_2, [0.3, 0.8], [0.3, 0.9], [0.1, 0.2], [-0.7, -0.6], [-0.6, -0.2], [-0.6, -0.2], 0.7\rangle, \\ \langle x_3, [0.4, 0.7], [0.5, 0.7], [0.3, 0.4], [-0.9, -0.5], [-0.4, -0.3], [-0.8, -0.1], 0.3\rangle \end{cases}$$

and

$$A_2 = \begin{cases} \langle x_1, [0.2,0.8], [0.3,0.6], [0.3,0.6], [-0.3,-0.2], [-0.7,-0.5], [-0.5,-0.4], 0.1 \rangle, \\ \langle x_2, [0.4,0.7], [0.5,0.7], [0.2,0.3], [-0.2,-0.1], [-0.8,-0.4], [-0.9,-0.8], 0.6 \rangle, \\ \langle x_3, [0.5,0.6], [0.3,0.5], [0.1,0.4], [-0.4,-0.2], [-0.8,-0.5], [-0.7,-0.6], 0.9 \rangle \end{cases}$$ are two IVBFWN-sets in $X$.

Then their union is given as follows:

$$A_1 \cup A_2 = \begin{cases} \langle x_1, [0.3,0.9], [0.2,0.7], [0.2,0.5], [-0.8,-0.7], [-0.6,-0.3], [-0.4,-0.3], 0.5 \rangle, \\ \langle x_2, [0.4,0.8], [0.4,0.8], [0.1,0.2], [-0.7,-0.6], [-0.7,-0.3], [-0.6,-0.2], 0.7 \rangle, \\ \langle x_3, [0.5,0.7], [0.4,0.6], [0.1,0.4], [-0.9,-0.5], [-0.6,-0.4], [-0.7,-0.1], 0.9 \rangle \end{cases}$$

Then their intersection is given as follows:

$$A_1 \cap A_2 = \begin{cases} \langle x_1, [0.2,0.8], [0.2,0.7], [0.3,0.6], [-0.3,-0.2], [-0.6,-0.3], [-0.5,-0.4], 0.1 \rangle, \\ \langle x_2, [0.3,0.7], [0.4,0.8], [0.2,0.3], [-0.2,-0.1], [-0.7,-0.3], [-0.9,-0.8], 0.6 \rangle, \\ \langle x_3, [0.4,0.6], [0.4,0.6], [0.3,0.4], [-0.4,-0.2], [-0.6,-0.4], [-0.8,-0.6], 0.3 \rangle \end{cases}$$

Then the complement of $A_1$ is given as follows:

$$A_1^C = \begin{cases} \langle x_1, [0.2,0.5], [0.2,0.9], [0.3,0.9], [-0.4,-0.3], [-0.9,-0.5], [-0.8,-0.7], 0.5 \rangle, \\ \langle x_2, [0.1,0.2], [0.1,0.7], [0.3,0.8], [-0.6,-0.2], [-0.7,-0.3], [-0.7,-0.6], 0.3 \rangle, \\ \langle x_3, [0.3,0.4], [0.3,0.5], [0.4,0.7], [-0.8,-0.1], [-0.7,-0.6], [-0.9,-0.5], 0.7 \rangle \end{cases}$$

Note that the a IVBFWN-number is denoted
$$\tilde{a} = \langle [T_L^+, T_R^+], [I_L^+, I_R^+], [F_L^+, F_R^+], [T_L^-, T_R^-],$$
$$[I_L^-, I_R^-], [F_L^-, F_R^-], p \rangle$$
for convenience.

**Definition 3.6.** Let
$$A_1 = \langle [T_{1L}^+, T_{1R}^+], [I_{1L}^+, I_{1R}^+], [F_{1L}^+, F_{1R}^+],$$
$$[T_{1L}^-, T_{1R}^-], [I_{1L}^-, I_{1R}^-], [F_{1L}^-, F_{1R}^-], p_1 \rangle$$
and
$$A_2 = \langle [T_{2L}^+, T_{2R}^+], [I_{2L}^+, I_{2R}^+], [F_{2L}^+, F_{2R}^+],$$
$$[T_{2L}^-, T_{2R}^-], [I_{2L}^-, I_{2R}^-], [F_{2L}^-, F_{2R}^-], p_2 \rangle$$
be two IVBFWN-number. Then the operations for IVBFWN-numbers are defined as below;

i. $\lambda A_1 = \langle [1-(1-T_L^+)^\lambda, 1-(1-T_R^+)^\lambda],$
$$[(I_L^+)^\lambda, (I_R^+)^\lambda], [(F_L^+)^\lambda, (F_R^+)^\lambda],$$
$$[-(-T_L^-)^\lambda, -(-T_R^-)^\lambda], [-(-I_L^-)^\lambda, -(-I_R^-)^\lambda],$$
$$[-(1-(1-(-F_L^-))^\lambda), -(1-(1-(-F_R^-))^\lambda)], p_1 \rangle$$

ii. $A_1^\lambda = \langle [(T_L^+)^\lambda, (T_R^+)^\lambda], [1-(1-I_L^+)^\lambda, 1-(1-I_R^+)^\lambda],$
$$[1-(1-F_L^+)^\lambda, 1-(1-F_R^+)^\lambda], [-(1-(1-(-T_L^-))^\lambda),$$
$$-(1-(1-(-T_R^-))^\lambda)], [-(-I_L^-)^\lambda, -(-I_R^-)^\lambda],$$
$$[-(-F_L^-)^\lambda, -(-F_R^-)^\lambda], p_1 \rangle$$

iii. $A_1 + A_2 = \langle [T_{1L}^+ + T_{2L}^+ - T_{1L}^+ T_{2L}^+, T_{1R}^+ + T_{2R}^+ - T_{1R}^+ T_{2R}^+],$
$$[I_{1L}^+ I_{2L}^+, I_{1R}^+ I_{2R}^+], [F_{1L}^+ F_{2L}^+, F_{1R}^+ F_{2R}^+], [-T_{1L}^- T_{2L}^-,$$
$$-T_{1R}^- T_{2R}^-], [-(-I_{1L}^- - I_{2L}^- - I_{1L}^- I_{2L}^-),$$
$$-(-I_{1R}^- - I_{2R}^- - I_{1R}^- I_{2R}^-)], [-(-F_{1L}^- - F_{2L}^- - F_{1L}^- F_{2L}^-),$$
$$-(-F_{1R}^- - F_{2R}^- - F_{1R}^- F_{2R}^-)], \max\{p_1, p_2\} \rangle$$

iv. $A_1 A_2 = \langle [T_{1L}^+ T_{2L}^+, T_{1R}^+ T_{2R}^+],$
$$[I_{1L}^+ + I_{2L}^+ - I_{1L}^+ I_{2L}^+, I_{1R}^+ + I_{2R}^+ - I_{1R}^+ I_{2R}^+],$
$$[F_{1L}^+ + F_{2L}^+ - F_{1L}^+ F_{2L}^+, F_{1R}^+ + F_{2R}^+ - F_{1R}^+ F_{2R}^+],$
$$[-(-T_{1L}^- - T_{2L}^- - T_{1L}^- T_{2L}^-), -(-T_{1R}^- - T_{2R}^- - T_{1R}^- T_{2R}^-)],$
$$[I_{1L}^- I_{2L}^-, I_{1R}^- I_{2R}^-], [F_{1L}^- F_{2L}^-, F_{1R}^- F_{2R}^-], \min\{p_1, p_2\} \rangle$$

where $\lambda > 0$.

**Definition 3.7.** Let
$$\tilde{a} = \langle [T_L^+, T_R^+], [I_L^+, I_R^+], [F_L^+, F_R^+], [T_L^-, T_R^-],$$
$$[I_L^-, I_R^-], [F_L^-, F_R^-], p \rangle$$

be a IVBFWN-number. Then, the score function $S(\tilde{a})$ accuracy function $A(\tilde{a})$ and certainty function $C(\tilde{a})$ of an NBN are defined as follows:

$$S(\tilde{a}) = \frac{p}{12}\big(T_L^+ + T_R^+ + 1 - I_L^+ + 1 - I_R^+ + 1 - F_L^+ + 1 - F_R^+ +$$
$$1 + T_L^- + 1 + T_R^- - I_L^- - I_R^- - F_L^- - F_R^-\big)$$

$$A(\tilde{a}) = \frac{p}{8}(4 + T_L^+ + T_R^+ - F_L^+ - F_R^+ + T_L^- + T_R^- - F_L^- - F_R^-)$$

$$C(\tilde{a}) = \frac{p}{4}(2 + T_L^+ + T_R^+ - F_L^- - F_R^-)$$

The comparison method can be defined as follows:

i. If $S(\tilde{a}_1) > S(\tilde{a}_2)$, then $\tilde{a}_1$ is greater than $\tilde{a}_2$, that is, $\tilde{a}_1$ is superior to $\tilde{a}_2$, denoted by $\tilde{a}_1 > \tilde{a}_2$;

ii. If $S(\tilde{a}_1) = S(\tilde{a}_2)$, and $A(\tilde{a}_1) > A(\tilde{a}_2)$, then $\tilde{a}_1$ is greater than $\tilde{a}_2$, that is, $\tilde{a}_1$ is superior to $\tilde{a}_2$, denoted by $\tilde{a}_1 < \tilde{a}_2$;

iii. If $S(\tilde{a}_1) = S(\tilde{a}_2)$, $A(\tilde{a}_1) = A(\tilde{a}_2)$, and $C(\tilde{a}_1) > C(\tilde{a}_2)$, then $\tilde{a}_1$ is greater than $\tilde{a}_2$, that is, $\tilde{a}_1$ is superior to $\tilde{a}_2$, denoted by $\tilde{a}_1 > \tilde{a}_2$;

iv.  If $S(\tilde{a}_1) = S(\tilde{a}_2)$, $A(\tilde{a}_1) > A(\tilde{a}_2)$, and $C(\tilde{a}_1) = C(\tilde{a}_2)$, then $\tilde{a}_1$ is equal to $\tilde{a}_2$, that is, $\tilde{a}_1$ is indifferent to $\tilde{a}_2$, denoted by $\tilde{a}_1 > \tilde{a}_2$;

**Definition 3.8.** Let

$\tilde{a}_j = \langle [T_{jL}^+, T_{jR}^+], [I_{jL}^+, I_{jR}^+], [F_{jL}^+, F_{jR}^+], [T_{jL}^-, T_{jR}^-],$
$[I_{jL}^-, I_{jR}^-], [F_{jL}^-, F_{jR}^-], p_j \rangle (j=1,2,...,n)$

be a family of IVBFWN-numbers. A mapping $A_p : \Im_n \to \Im$ is called IVBFWN average operator if it satisfies

$A_p(\tilde{a}_1, \tilde{a}_2, ..., \tilde{a}_n) = \sum_{j=1}^{n} p_j \tilde{a}_j$

$= \left\langle \left[ 1 - \prod_{j=1}^{n}(1-T_{jL}^+)^{p_j}, 1 - \prod_{j=1}^{n}(1-T_{jR}^+)^{p_j} \right], \right.$

$\left[ \prod_{j=1}^{n}(I_{jL}^+)^{p_j}, \prod_{j=1}^{n}(I_{jR}^+)^{p_j} \right], \left[ \prod_{j=1}^{n}(F_{jL}^+)^{p_j}, \prod_{j=1}^{n}(F_{jR}^+)^{p_j} \right],$

$\left[ -\prod_{j=1}^{n}(-T_{jL}^-)^{p_j}, -\prod_{j=1}^{n}(-T_{jR}^-)^{p_j} \right],$

$\left[ -\left(1 - \prod_{j=1}^{n}(1-(-I_{jL}^-))^{p_j}\right), -\left(1 - \prod_{j=1}^{n}(1-(-I_{jR}^-))^{p_j}\right) \right],$

$\left. \left[ -\left(1 - \prod_{j=1}^{n}(1-(-F_{jL}^-))^{p_j}\right), -\left(1 - \prod_{j=1}^{n}(1-(-F_{jR}^-))^{p_j}\right) \right], \max_j\{p_j\} \right\rangle$

**Theorem 3.9.** Let

$\tilde{a}_j = \langle [T_{jL}^+, T_{jR}^+], [I_{jL}^+, I_{jR}^+], [F_{jL}^+, F_{jR}^+], [T_{jL}^-, T_{jR}^-],$
$[I_{jL}^-, I_{jR}^-], [F_{jL}^-, F_{jR}^-], p_j \rangle (j=1,2,...,n)$

be a family of IVBFWN-numbers. Then,
i.  If $\tilde{a}_j = \tilde{a}$ for all $j=1,2,...,n$ then, $A_p(\tilde{a}_1, \tilde{a}_2, ..., \tilde{a}_n) = \tilde{a}$
ii.  $\min_{j=1,2,...,n} \tilde{a}_j \leq A_p(\tilde{a}_1, \tilde{a}_2, ..., \tilde{a}_n) \leq \max_{j=1,2,...,n} \tilde{a}_j$
iii. If $\tilde{a}_j = \tilde{a}_j^*$ for all $j=1,2,...,n$ then,
$A_p(\tilde{a}_1, \tilde{a}_2, ..., \tilde{a}_n) \leq A_p(\tilde{a}_1^*, \tilde{a}_2^*, ..., \tilde{a}_n^*)$

**Definition 3.10.** Let

$\tilde{a}_j = \langle [T_{jL}^+, T_{jR}^+], [I_{jL}^+, I_{jR}^+], [F_{jL}^+, F_{jR}^+], [T_{jL}^-, T_{jR}^-],$
$[I_{jL}^-, I_{jR}^-], [F_{jL}^-, F_{jR}^-], p_j \rangle (j=1,2,...,n)$

be a family of IVBFWN-numbers. A mapping $G_p : \Im_n \to \Im$ is called IVBFWN geometric operator if it satisfies

$G_p(\tilde{a}_1, \tilde{a}_2, ..., \tilde{a}_n) = \prod_{j=1}^{n} \tilde{a}_j^{p_j}$

$= \left\langle \left[ \prod_{j=1}^{n}(T_{jL}^+)^{p_j}, \prod_{j=1}^{n}(T_{jR}^+)^{p_j} \right], \right.$

$\left[ 1 - \prod_{j=1}^{n}(1-I_{jL}^+)^{p_j}, 1 - \prod_{j=1}^{n}(1-I_{jR}^+)^{p_j} \right],$

$\left[ 1 - \prod_{j=1}^{n}(1-F_{jL}^+)^{p_j}, 1 - \prod_{j=1}^{n}(1-F_{jR}^+)^{p_j} \right],$

$\left[ -\left(1 - \prod_{j=1}^{n}(1-(-T_{jL}^-))^{p_j}\right), -\left(1 - \prod_{j=1}^{n}(1-(-T_{jR}^-))^{p_j}\right) \right],$

$\left[ -\prod_{j=1}^{n}(-I_{jL}^-)^{p_j}, -\prod_{j=1}^{n}(-I_{jR}^-)^{p_j} \right],$

$\left. \left[ -\prod_{j=1}^{n}(-F_{jL}^-)^{p_j}, -\prod_{j=1}^{n}(-F_{jR}^-)^{p_j} \right], \min_j\{p_j\} \right\rangle$

**Theorem 3.11.** Let

$\tilde{a}_j = \langle [T_{jL}^+, T_{jR}^+], [I_{jL}^+, I_{jR}^+], [F_{jL}^+, F_{jR}^+], [T_{jL}^-, T_{jR}^-],$
$[I_{jL}^-, I_{jR}^-], [F_{jL}^-, F_{jR}^-], p_j \rangle (j=1,2,...,n)$

be a family of IVBFWN-numbers. Then,
i.  If $\tilde{a}_j = \tilde{a}$ for all $j=1,2,...,n$ then, $G_p(\tilde{a}_1, \tilde{a}_2, ..., \tilde{a}_n) = \tilde{a}$
ii.  $\min_{j=1,2,...,n} \tilde{a}_j \leq G_p(\tilde{a}_1, \tilde{a}_2, ..., \tilde{a}_n) \leq \max_{j=1,2,...,n} \tilde{a}_j$
iii. If $\tilde{a}_j = \tilde{a}_j^*$ for all $j=1,2,...,n$ then,
$G_P(\tilde{a}_1, \tilde{a}_2, ..., \tilde{a}_n) \leq G_P(\tilde{a}_1^*, \tilde{a}_2^*, ..., \tilde{a}_n^*)$

Note that the aggregation results are still NBNs

## IV. NBN- DECISION MAKING METHOD

In this section, we develop an approach based on the $A_P$ (or $G_P$) operator and the above ranking method to deal with multiple criteria decision making problems with IVBFWN-information.

Suppose that $A = \{A_1, A_2, ..., A_m\}$ and $C = \{C_1, C_2, ..., C_n\}$ is the set of alternatives and criterions or attributes, respectively. Let $p_j$ be the fuzzy weight index of attributes, such that $p_j \in [0,1](j=1,2,...,n)$ and $p_j$ refers to the weight of attribute $C_j$. An alternative on criterions is evaluated by the

decision maker, and the evaluation values are represented by the form of IVBFWN-numbers. Assume that

$$(\tilde{a}_{ij})_{m \times n} = (\langle [T_{ijL}^+, T_{ijR}^+], [I_{ijL}^+, I_{ijR}^+], [F_{ijL}^+, F_{ijR}^+], [T_{ijL}^-, T_{ijR}^-],$$
$$[I_{ijL}^-, I_{ijR}^-], [F_{ijL}^-, F_{ijR}^-], p_{ij} \rangle)_{m \times n}$$

is the decision matrix provided by the decision maker; $\tilde{a}_{ij}$ is a IVBFWN-number for alternative $A_i$ associated with the criterions $C_j$. We have the conditions

$$T_{ijL}^+, T_{ijR}^+, I_{ijL}^+, I_{ijR}^+, F_{ijL}^+, F_{ijR}^+, T_{ijL}^-, T_{ijR}^-, I_{ijL}^-, I_{ijR}^-, F_{ijL}^-, F_{ijR}^- \in [0,1]$$

such that

$$0 \leq T_{ijL}^+ + T_{ijR}^+ + I_{ijL}^+ + I_{ijR}^+ + F_{ijL}^+ + F_{ijR}^+ -$$
$$T_{ijL}^- - T_{ijR}^- - I_{ijL}^- - I_{ijR}^- - F_{ijL}^- - F_{ijR}^- \leq 12$$

for $(i = 1, 2, ..., m)$ and $(j = 1, 2, ..., n)$.

Now, we can develop an algorithm as follows;

**Algorithm**

**Step 1.** Construct the decision matrix provided by the decision maker as;

$$(\tilde{a}_{ij})_{m \times n} = (\langle [T_{ijL}^+, T_{ijR}^+], [I_{ijL}^+, I_{ijR}^+], [F_{ijL}^+, F_{ijR}^+], [T_{ijL}^-, T_{ijR}^-],$$
$$[I_{ijL}^-, I_{ijR}^-], [F_{ijL}^-, F_{ijR}^-], p_{ij} \rangle)_{m \times n}$$

**Step 2.** Compute $\tilde{a}_i = A_p(\tilde{a}_{i1}, \tilde{a}_{i2}, ..., \tilde{a}_{in})$ (or $G_p(\tilde{a}_{i1}, \tilde{a}_{i2}, ..., \tilde{a}_{in})$) for each $\tilde{a}_i (i = 1, 2, ..., m)$

**Step 3.** Calculate the score values of $S(\tilde{a}_i)$ for the $(i = 1, 2, ..., m)$ collective overall IVBFWN-number of $\tilde{a}_i (i = 1, 2, ..., m)$

**Step 4.** Rank all the software systems of $\tilde{a}_i (i = 1, 2, ..., m)$ according to the score values

Now, we give a numerical example as follows;

**Example 4.1.** Let us consider decision making problem adapted from Ye [16]. There is an investment company, which wants to invest a sum of money in the best option. There is a panel with the set of the four alternatives is denoted by $C_1$ = car company $C_2$ = food company, $C_3$ = computer company, $C_4$ = arms company to invest the money. The investment company must take a decision according to the set of the four attributes is denoted by $A_1$ = risk, $A_2$ = growth, $A_3$ = environmental impact, $A_4$ = performance. Then the according to this algorithm, we have,

**Step 1.** Construct the decision matrix provided by the customer as;

**Table 1:** Decision matrix given by customer

| | $C_1$ |
|---|---|
| $A_1$ | $\langle [0.5, 0.6], [0.2, 0.5], [0.1, 0.7], [-0.2, -0.1], [-0.6, -0.2], [-0.4, -0.3], 0.5 \rangle$ |
| $A_2$ | $\langle [0.1, 0.2], [0.3, 0.8], [0.2, 0.4], [-0.5, -0.2], [-0.9, -0.3], [-0.6, -0.1], 0.8 \rangle$ |
| $A_3$ | $\langle [0.4, 0.8], [0.4, 0.6], [0.4, 0.6], [-0.3 - 0.2], [-0.7, -0.5], [-0.5, -0.4], 0.2 \rangle$ |
| $A_4$ | $\langle [0.6, 0.9], [0.3, 0.8], [0.5, 0.6], [-0.8, -0.5], [-0.5, -0.1], [-0.2, -0.1], 0.3 \rangle$ |

| | $C_2$ |
|---|---|
| $A_1$ | $\langle [0.3, 0.9], [0.1, 0.8], [0.2, 0.5], [-0.8, -0.7], [-0.5, -0.1], [-0.4, -0.1], 0.6 \rangle$ |
| $A_2$ | $\langle [0.2, 0.8], [0.1, 0.4], [0.3, 0.4], [-0.5, -0.1], [-0.3, -0.1], [-0.9, -0.4], 0.4 \rangle$ |
| $A_3$ | $\langle [0.1, 0.6], [0.3, 0.9], [0.3, 0.5], [-0.8, -0.7], [-0.4, -0.3], [-0.7, -0.6], 0.7 \rangle$ |
| $A_4$ | $\langle [0.1, 0.2], [0.8, 0.9], [0.2, 0.7], [-0.5, -0.4], [-0.6, -0.3], [-0.5, -0.3], 0.1 \rangle$ |

| | $C_3$ |
|---|---|
| $A_1$ | $\langle [0.1, 0.6], [0.1, 0.5], [0.1, 0.4], [-0.5, -0.2], [-0.7, -0.3], [-0.4, -0.2], 0.9 \rangle$ |
| $A_2$ | $\langle [0.3, 0.4], [0.1, 0.6], [0.5, 0.7], [-0.5, -0.1], [-0.8, -0.7], [-0.9, -0.8], 0.3 \rangle$ |
| $A_3$ | $\langle [0.3, 0.9], [0.2, 0.8], [0.2, 0.3], [-0.5, -0.4], [-0.6, -0.5], [-0.7, -0.6], 0.5 \rangle$ |
| $A_4$ | $\langle [0.2, 0.7], [0.5, 0.8], [0.8, 0.9], [-0.9, -0.8], [-0.8, -0.5], [-0.5, -0.2], 0.4 \rangle$ |

| | $C_4$ |
|---|---|
| $A_1$ | $\langle [0.6, 0.8], [0.4, 0.6], [0.1, 0.3], [-0.4, -0.3], [-0.6, -0.3], [-0.7, -0.5], 0.7 \rangle$ |
| $A_2$ | $\langle [0.3, 0.8], [0.3, 0.9], [0.1, 0.2], [-0.8, -0.6], [-0.6, -0.4], [-0.4, -0.2], 0.1 \rangle$ |
| $A_3$ | $\langle [0.7, 0.9], [0.1, 0.4], [0.2, 0.6], [-0.7, -0.6], [-0.9, -0.5], [-0.3, -0.2], 0.2 \rangle$ |
| $A_4$ | $\langle [0.4, 0.6], [0.3, 0.5], [0.1, 0.7], [-0.3, -0.1], [-0.6, -0.5], [-0.7, -0.3], 0.8 \rangle$ |

**Step 2.** Compute $\tilde{a}_i = A_p(\tilde{a}_{i1}, \tilde{a}_{i2}, \tilde{a}_{i3}, \tilde{a}_{i4})$ for each $(i = 1, 2, 3, 4)$ as;

$\tilde{a}_1$  $\langle [0.4, 0.8], [0.2, 0.6], [0.1, 0.5], [-0.4, -0.3], [-0.6, -0.2], [-0.5, -0.3], 0.9 \rangle$

$\tilde{a}_2$  $\langle [0.2, 0.6], [0.2, 0.6], [0.2, 0.4], [-0.6, -0.2], [-0.7, -0.4], [-0.8, -0.5], 0.8 \rangle$

$\tilde{a}_3$  $\langle [0.4, 0.8], [0.2, 0.7], [0.3, 0.5], [-0.5, -0.4], [-0.7, -0.5], [-0.6, -0.5], 0.7 \rangle$

$\tilde{a}_4$  $\langle [0.3, 0.7], [0.4, 0.7], [0.3, 0.7], [-0.6, -0.4], [-0.6, -0.4], [-0.5, -0.2], 0.8 \rangle$

**Step 3.** Calculate the score values of $S(\tilde{a}_i)(i=1,2,3,4)$ for the collective overall IVBFWN-number of $\tilde{a}_i (i=1,2,...,m)$ as;

$S(\tilde{a}_1) = 0.50 \quad S(\tilde{a}_2) = 0.47 \quad S(\tilde{a}_3) = 0.40 \quad S(\tilde{a}_4) = 0.37$

**Step 4.** Rank all the software systems of $A_i (i=1,2,3,4)$ according to the score values as;
$A_1 \succ A_2 \succ A_3 \succ A_4$
and thus $A_1$ is the most desirable alternative.

## Conclusion

This paper presented an interval-valued bipolar neutrosophic set and its score, certainty and accuracy functions. In the future, we shall further study more aggregation operators for interval-valued bipolar neutrosophic set and apply them to solve practical applications in group decision making, expert system, information fusion system, game theory, and so on.